\documentclass[10pt, a4paper]{amsart}
\input amssym.def
\input amssym.tex
\title[ Ranks
 of the Sylow 2-Subgroups of the 
 Classical Groups     ]
{ $\mbox { Ranks
 of the Sylow 2-Subgroups of the 
 Classical Groups    }$   }   

\author{ Mong Lung Lang}

\begin{document}
\baselineskip=11pt


\maketitle

\vspace {-1cm}

\begin{abstract}
Let $S$ be a  2-group. The rank (normal rank)
of $S$ is the maximal dimension of an elementary
 abelian subgroup (a normal
 elementary
 abelian subgroup) of $S$ over $\Bbb Z_2$.
 The purpose of this article is to determine
 the  rank and normal rank of $S$,
 where $S$ is a Sylow 2-subgroup of the 
 classical groups of odd characteristic.
\end{abstract} 

 \section{ Introduction} 
 Sylow 2-subgroups of
 the classical groups 
 (of odd characteristic) have been
 completely determined by Carter and Fong [CF] and
 Wong [W].
Let $S$ be such a 2-group. The rank (normal rank)
of $S$ is the maximal dimension of an elementary
 abelian subgroup (a normal
 elementary
 abelian subgroup) of $S$ over $\Bbb Z_2$.
 The purpose of this article is to determine
 the  rank and normal rank of $S$.
 The result of this article
 is tabulated in the following table.
 The remaining of this article is organised 
 as follows :

\smallskip
Section 2 gives some very basic facts about 
 ranks and normal ranks of direct product and 
 semidirect product of groups. Section 3 studies the 
 ranks and normal ranks of wreath
 products. These results enable us to determine
 the ranks and normal ranks of Sylow 2-subgroups 
 of  $GL_n(q)$, $U_n(q)$, $O^+_{2n+1}
(q)$, $O_{2n}(\eta, q)$, $SL_{2n+1}(q),$  
 $Sp_{2n}(q)$,  and $\Omega_{2n}(\eta, q)\,\,
(\eta = \pm 1, q^n \equiv -\eta)\,\,(\mbox{mod }4)$.
 An alternative  proof (cohomology free) of a lemma
 (Lemma 10.32 of [GLS2])
 which is very useful in the study of $p$-groups 
 is also provided.

\smallskip
Section 4 studies twisted wreath product of groups.
 Let $S(T,R,J)$ be a twisted wreath product
 (see section 4.2).
 We are able to give upper
 and lower bounds of $r_2(S(T,R,J))$.
 These  bounds are by no means the best possible
 bounds for arbitrary $T$'s and $R$'s
but will be proven in the section 5
 that are optimal if $S$ and $T$ are
 the ones associated to 
the Sylow 2-subgroups of the classical groups
 (of odd characteristic).
  Our study of 
 the twisted wreath product
 is motivated by the fact that   if
 $S$ is a Sylow 2-subgroup of the classical
 groups which cannot be described as a
 direct product of wreath products
 then
$S$ can be described by the usage of the 
twisted wreath product
 $S(T,R,J)$.

\smallskip Section 5 gives the ranks
 and normal ranks of the Sylow 2-subgroups
 of classical groups 
 that are described by twisted wreath
 product.  They are easy consequences
 of Propositions 4.4 and 4.9 except
 for $\Omega_{2(2m+1)}(\eta, q)$
 ($q^{2m+1} \equiv \eta$ (mod 4))
which requires some special treatment.
 Our analysis shows that 
 certain invariance of  $TR$ determines  uniquely
 the rank and normal rank of  $S(T,R,J)$.

\smallskip

The notations we used are basically those 
 in [CF] and [W].
In the following table, 
 $\eta = \pm 1$,
$\mbox{ord}_2(q^2-1)$
is the largest integer $m$ such that 
 $2^m$ is a divisor of $q^2-1$
 and that $[x]$ is the largest integer
 less than or equal to $x$.

\medskip
{\small \begin{center}
$
\begin{array}{|l|l|l|} \hline
 \vrule height 8pt width 0pt depth 5pt
\phantom{\Big |}
 \mbox{Group  }  & \mbox{2 Rank }  &
\mbox{Normal 2 rank}  \\ \hline
\phantom{\Big |}
SL_{2n} (q)\, :\, q \equiv 3 \,(4) & 2n-1 & n\\ \hline
\phantom{\Big |}
SL_{2n}(q)\, :\, q \equiv  1\,(4) & 2n-1 & 2n-1\\ \hline
\phantom{\Big |}
SU_{2n}(q)\, :\, q \equiv 1 \,(4) & 2n-1 & n\\ \hline
\phantom{\Big |}
SU_{2n}(q)\, :\, q \equiv  3\,(4) & 2n-1 & 2n-1\\\hline
\phantom{\Big |}
Sp_{2n}(q) & n & n\\
\hline
\phantom{\Big |}
\Omega _{2n+1}(q)\,:\
\mbox{ord}_2(q^2-1) \ge 4 & 2n & n\\
\hline
\phantom{\Big |}
\Omega _{2n+1}(q)\,:\
\mbox{ord}_2(q^2-1) = 3 & 2n & 2n\\
\hline
\phantom{\Big |}
\Omega _{2n}(\eta, q)\,:\
 q^n \equiv -\eta\,(4)\,\,\,
\mbox{ord}_2(q^2-1) \ge 4
 & 2n-2 & n-1\\
\hline
\phantom{\Big |}
\Omega _{2n}(\eta, q)\,:\
 q^n \equiv -\eta\,(4)\,\,\,
\mbox{ord}_2(q^2-1) =3
 & 2n-2 & 2n-2\\ \hline
\phantom{\Big |}
\Omega _{4n}(\eta, q)\,:\
 q^n \equiv \eta\,(4)\,\,\,
\mbox{ord}_2(q^2-1) \ge 4
 & 4n-1 & 2n\\
\hline
\phantom{\Big |}
\Omega _{4n}(\eta, q)\,:\
 q^n \equiv \eta\,(4)\,\,\,
\mbox{ord}_2(q^2-1) =3
 & 4n-1 & 4n-1\\
\hline

\phantom{\Big |}
\Omega _{4n+2}(\eta, q)\,:\
 q^n \equiv \eta\,(4)\,\,\,
\mbox{ord}_2(q^2-1) \ge 4
 & 4n+1 & 2n+1\\
\hline
\phantom{\Big |}
\Omega _{4n+2}(\eta, q)\,:\
 q^n \equiv \eta\,(4)\,\,\,
\mbox{ord}_2(q^2-1) =3
 & 4n+1 & 4n+1\\
\hline

\phantom{\Big |}
GL_n (q)\, :\, q \equiv 3 \,(4) & n & [(n+1)/2]\\ \hline
\phantom{\Big |}
GL_{n}(q)\, :\, q \equiv  1\,(4) & n & n\\ \hline
\phantom{\Big |}
U_{n}(q)\, :\, q \equiv 1 \,(4) & n &  [(n+1)/2]\\ \hline
\phantom{\Big |}
U_{n}(q)\, :\, q \equiv  3\,(4) & n & n\\\hline
\phantom{\Big |}
O_{2n+1}^+(q)\, :\,
\mbox{ord}_2(q^2-1) =3
  & 2n & 2n\\ \hline
\phantom{\Big |}
O_{2n+1}^+(q)\, :\,
\mbox{ord}_2(q^2-1) \ge 4
  & 2n & n\\ \hline
\phantom{\Big |}
O_{2n} (\eta, q)\,:\
 q^n \equiv \eta\,(4)\,\,\,
\mbox{ord}_2(q^2-1) =3
 & 2n & 2n\\
\hline
\phantom{\Big |}
O_{2n} (\eta, q)\,:\
 q^n \equiv \eta\,(4)\,\,\,
\mbox{ord}_2(q^2-1) \ge 4
 & 2n & n\\
\hline
\phantom{\Big |}
O_{2n} (\eta, q)\,:\
 q^n \equiv -\eta\,(4)\,\,\,
\mbox{ord}_2(q^2-1) =3
 & 2n & 2n\\
\hline
\phantom{\Big |}
O_{2n} (\eta, q)\,:\
 q^n \equiv -\eta\,(4)\,\,\,
\mbox{ord}_2(q^2-1) \ge 4
 & 2n & n+1\\
\hline

\end{array}$
\end{center}}

\section{Preliminaries : Ranks and
 Normal ranks of 2 groups}
 The main purpose of this section is to give some basic
 results that enable us to describe the ranks and normal
 ranks of 2-groups.

\smallskip

\noindent {\bf Definition 2.1.} 
Lat $p$ be a prime and let $G$ be a group.
 We say $G$ is of $p$-rank $n$ if  the $p$-rank of 
 a Sylow $p$-subgroup of $G$ is $n$.
 We say $G$ is of normal $p$-rank $n$ if    the 
normal rank of a Sylow $p$-subgroup of $G$ is $n$.
 The rank and normal rank of $G$ are  denoted by
 $r_p(G)$ and $nr_p(G)$ respectively.
Recall that 
 the rank (normal rank)
of a $p$-group  is the maximal dimension of an elementary
 abelian subgroup (a normal
 elementary
 abelian subgroup)  over $\Bbb Z_p$.

\smallskip
\noindent {\bf Lemma 2.2.} {\em Let $G= A \times B$. Then
 $r_p(G) = r_p(A) +r_p(B)$ and  $nr_p(G) = nr_p(A) +nr_p(B)$.}
\smallskip

\noindent {\em Proof.} Let  $E\subseteq G$
 be elementary abelian of  rank $r_p(G)$. 
 Define the projection of $E$ on $A$ and $B$ as follows.
 $$E|A = \{ a \in A\,:\, ab  \in E\mbox{ for some } b \in B\},
$$$$
E|B = \{ b \in B \,:\, ab  \in E\mbox{ for some } a \in A\}.$$
One sees easily that $E|A$ and $E|B$
 are elementary abelian subgroups of $A$ and $B$ respectively.
 Further,
 $$E \subseteq E|A \times E|B.\eqno (2.1)$$
 It follows easily from (2.1) that 
 $r_p(G) = r_p(A) +r_p(B)$. Note
 that if $E$ is normal, then both 
 $E|A$ and $E|B$ are normal in $A$ and $B$
 respectively as well. Hence
 $nr_p(G) = nr_p(A) +nr_p(B)$. \qed

\smallskip
\noindent {\bf Lemma 2.3.} {\em Let $S = A\rtimes B$ be the 
 semidirect product of $A$ and $B$
 and let $E$ be an  elementary abelian 
 subgroup of $S$. Suppose that 
 $E = (E \cap A)  \times R$. Then
 $r_2(E\cap A) \le r_2(A)$, $r_2(R) \le r_2(B)$,
$r_2(R) = r_2(R|B)$.
 In particular, 
 $r_2(S) \le  r_2(A) +r_2(B).$
 In the case $E$ is normal, the following holds.
$r_2(E\cap A) \le nr_2(A)$, $r_2(R) \le r_2(B)$,
 and 
 $nr_2(S) \le  nr_2(A) +r_2(B).$
}

\smallskip
\noindent {\em Proof.} Let $E$ be an elementary
 abelian 2-group of $S$.
Set $E = (E\cap A) \times R$.
Since $E\cap A$ is  in $A$, it is clear
 that 
$$r_2(E\cap A) \le r_2(A).\eqno (2.2)$$
 Note that 
 elements in $R^{\times} $ are of the forms $xb$, where
  $ x\in A, b\in B-\{1\}$.
Let 
$ M= \{x_1b_1, \cdots , x_nb_n\}$ be  a 
 minimal generating set of $R$.
Since elements in $M$ are of order 2,
$o(b_i) = 2$ for all $i$.
Since $x_ib_ix_jb_j =x_ix_j'b_ib_j $ is of order 2 as well,
 $o(b_ib_j) = 1$ or 2.
 This implies that 
 $ V = \left < b_1, b_2, \cdots , b_n\right
 >$ is an elementary abelian 2-subgroup of $B$. Further
 $R|B = V$.
 Suppose that 
$\{ b_1, b_2, \cdots , b_n\}$ is not a 
 minimal generating set of $V$.
 Without loss generality, we may assume that 
 $b_1 = \prod_{i=2}^n b_i^{e_i}$, where $e_i = 1$ or 2.
It follows that 
$$\prod_{i=2}^n (x_ib_i)^{e_i}
 = x_1'b_1.$$
Since $M$ is a minimal generating set of $R$,
 $\left <x_1'b_1, x_1b_1\right >
 = \left <x_1'b_1\right >
 \times \left < x_1b_1\right > \subseteq
 R$ is elementary
 abelian of order 4.
 Note that 
 $$ \left <x_1'b_1\right >
 \times \left < x_1b_1\right >
 = \left <x_1'b_1\right >\times 
\left < x_1b_1x_1'b_1   \right >\subseteq R .$$
 As the group $  \left <x_1b_1x_1'b_1\right >$  
 (of order 2) is in $E\cap A$,
 we have just shown that $R$ intersects $E\cap A$
 nontrivially. A contradiction.
 Hence $\{ b_1, b_2, \cdots , b_n\}$
 is a minimal generating set
 of $V \subseteq B$.
It follows that  $r_2(R) = |M|= n= r_2(V) = r_2(R|B)$
 and that 
  $r_2(R) = n \le r_2(B)$.
 As a consequence,
 $r_2(E) = r_2(E\cap A) + r_2(R)
 \le r_2(A) + r_2(B)$.

 In the case $E$ is normal, 
$E\cap A$ is a normal subgroup of $A$. Hence
$r_2(E\cap A) \le nr_2(A).$ 
This 
 completes the proof  of the lemma.
 \qed

\section {Wreath Product}

We shall first give a (cohomology free) proof of the following
 fact which is very useful in our study of the 
normal ranks of $p$-groups.
The purpose of such  proof is to keep our study of 
 the 2-ranks elementary.
 Our main results 
 concerning the ranks and normal ranks can be 
 found in Proposition 3.4.

\smallskip
\noindent {\bf Lemma 3.1.} 
 (Lemma 10.32 of [GLS2])
{\em Let  $P = Q_1 \times Q_2 \times 
 \cdots \times Q_p$, where the $Q_i$ are cycled by the element 
 $x$ of order $p$ $($whence $P \left <x\right
 > \cong  Q_1 \wr \Bbb Z_p$, the wreath 
 product of $Q_1$ and $\Bbb Z_p)$.
 Then every element of order $p$ in the coset $xP$ 
 is a $P$-conjugate of $x$. }
\smallskip

\noindent {\em Proof.}
We note first that every $P$-conjugate
 $pxp^{-1} $ are members in $xP$. 
This follows from the fact that $P$ is 
 a normal subgroup of $P\left <x \right >$.

Note also that $q_1q_2 \cdots q_p \in C_P(x)$
 $(q_i \in Q_i$)
 if and only if 
$xq_1q_2 \cdots q_px^{-1} = q_1q_2 \cdots q_p.$
 As $P$ is the 
 direct product of the $Q_i$'s, we must have
$$xq_px^{-1} = q_{1}, 
xq_ix^{-1} = q_{i+1},\mbox{ where }  i= 1,2,\cdots , p-1.$$
 As a consequence, there are $|Q_1|$ choices 
 of $q_1$ and the remaining $q_i$'s $(i \ge 2)$ are 
 uniquely determined by $q_1$ and $x$.
 Hence 
$$ |C_P(x)| = |Q_1|.$$
 This implies that $x$ possesses 
$|P|/ |C_P(x)| = |Q_1|^{p-1}$
$P$-conjugates in $xP$. Note that 
 such conjugates are elements of order $p$.

In order to complete our lemma, it suffices to prove
 that $Px$ has  $|Q_1|^{p-1}$ elements of order $p$.
Note that $q_1q_2 \cdots q_p x$ is of order $p$
 if and only if  $(q_1q_2 \cdots q_p x)^p =1$.
 It follows that 
$$(q_1q_2 \cdots q_p)(q_1q_2 \cdots q_p)^x
(q_1q_2 \cdots q_p)^{x^2}\cdots 
(q_1q_2 \cdots q_p)^{x^{p-1}}=1,$$
 where $q^x = xqx^{-1}$,
$q_1, q_p^x, q_{p-1}^{x^2},\cdots , q_2^{x^{p-1}}\in Q_1$.
Since the product $\prod Q_i$ is a direct product,
 we must have
 $$q_1q_p^x
 q_{p-1} ^{x^2}\cdots 
q_2 ^{x^{p-1}}=1.$$
 It follows that the choices of $q_1$ is uniquely
 determines by $x$ and $ q_2, q_3, \cdots , q_{p-1}$.
 Hence $xP$ possesses at most $|Q_1|^{p-1}$
 elements of order $p$.
 This completes the proof of our lemma.\qed

\noindent {\bf Remark.}
(i) Denote by $d(X)$ the number of elements of order $p$ of 
 $X$. Our lemma implies that $x^iP$ $(1\le
 i \le p-1)$ possesses 
 $|Q_1|^{p-1}$ elements of order $p$. 
 It is clear that $P$ possesses  $(d(Q_1)+1)^p -1 $ elements 
 of order $p$. As a consequence,
$$d(P\left<x\right >) = (d(Q_1)+1)^p -1 +  (p-1)|Q_1|^{p-1}.$$
 (ii) In the case 
 $Q_1$ is a Sylow $p$-subgroup of  $S_{p^n}$,
 the cycle decomposition of $x$ is $p^{n+1}$ and 
 $P\left <x \right >$ is a 
 Sylow $p$-subgroup of $S_{p^{n+1}}$.
Denote by $v_{n+1}$ the number of elements
 of  $P\left <x \right >$ with cycle decomposition
 $p^{n+1}$. Our lemma implies that $x^iP$ $(1\le
 i \le p-1)$ possesses 
 $|Q_1|^{p-1}$ elements with cycle decomposition $p^{n+1}$.
 It is clear that $P$ possesses  $v_n^p $ elements 
 with cycle decomposition $p^{n+1}$.
As a consequence, we have
$$v_{n+1} = v_n^p  + (p-1)|Q_1|^{p-1}.\eqno(3.1)$$
\smallskip

\noindent {\bf Lemma 3.2.}
 {\em  Let  $S = (Q_1\times Q_2)\rtimes \left <x\right >
 \cong Q_1 \wr \Bbb Z_2$
 be a $2$-group.
 Then $r_2(S) = 2\cdot r_2(Q_1)$.
 Suppose that the rank of $Q_1$ is at least $2$.
 Then every elementary abelian subgroup of
  dimension   $r_2(S)$
  is contained in $Q_1\times Q_2$.
}

\smallskip
\noindent {\em Proof.} Let $E$ be an elementary
 abelian subgroup of $S$ of dimension  $r_2(S)$.
 Set $E = A \times B$, where $A = E\cap (Q_1\times Q_2)$.
Applying Lemma 2.3, $r_2(B) \le 1$.
 In the case $r_2(B) = 0$, $E$ is a subgroup of
 $Q_1\times Q_2$. The maximality of $E$
 (in dimension) implies that $r_2(E) = 2 \cdot r_2(Q_1)$.
 Hence $r_2(S) = r_2(E) = 2 \cdot r_2(Q_1)$.
 In the case $r_2(B) =1$, $B$ is generated
 by an element of the form $ y= q_1q_2 x$, where
 $q_i \in Q_i$.
Note that $Q_1 \cap Q_2=1$,
 $[Q_1,Q_2]=1$ and that the conjugation action
 of $y$ interchanges $Q_1$ and $Q_2$.
  Let $t_1t_2\in E \cap (Q_1\times Q_2)$,
 where $t_i \in Q_i$. Since $E$ is abelian,
 $$[y, t_1t_2]=1.\eqno(3.2)$$
 It follows that 
$$yt_1y^{-1}  = t_2\mbox{ and }
y^{-1} t_2 y = t_1. \eqno (3.3)$$
Hence $t_2$ (of order 2)
is uniquely determined by $t_1$ (of order 2) and 
 $ y$. As a consequence,
 $r_2(E\cap (Q_1\times Q_2)) \le
 r_2(Q_1).$
 It follows that 
 $$r_2(S) = r_2(E) \le 1+ r_2(Q_1).\eqno(3.4)$$
 Since $Q_1\cong Q_2$ ($x$ permutes $Q_1$ and $Q_2$),
 the rank of $Q_1\times Q_2$ is $2\cdot r_2(Q_1)$.
 This together with (3.4) gives
$$ 2\cdot r_2(Q_1) \le r_2(S)  \le  1+ r_2(Q_1).$$
Hence $$r_2(Q_1) = r_2(Q_2) = 1. \eqno(3.5)$$
It follows that 
 $r_2(S) = 2 = 2\cdot r_2(Q)$.
 This completes the proof of the first part
 of our lemma.

\smallskip
Suppose that $S$ possesses an elementary abelian subgroup
 $E$ such that $E$ is not a subgroup of $Q_1\times Q_2$.
 It follows from the above that $r_2(Q_1) = 1.$
 A contradiction. Hence  $E \subseteq Q_1\times Q_2.$\qed

\smallskip
\noindent 
 {\bf Remark.} Suppose that $E\subseteq S$ is elementary
 abelian of dimension  $2r_2(Q_1) -1$, where $r(Q_1)\ge 2$.
 Applying the proof of our lemma (3.2), one can show that
$E$ is a subgroup of $Q_1\times Q_2$ as well.

\smallskip
\noindent {\bf Lemma 3.3.} {\em Let $S$ be a $2$-group
 of the form $(Q_1\times Q_2)\rtimes
\left <x\right > \cong Q_1 \wr \Bbb Z_2$.
Then $nr_2(S) = 2\cdot nr_2(Q_1)$. Further,
 the following hold.
\begin{enumerate}
\item[(i)]
Suppose that $Q_1$ is not elementary abelian.
 Then  every normal elementary abelian subgroup 
$($not necessarily of maximal dimension$)$ of $S$
  is contained in $Q_1\times Q_2$.
\item[(ii)] Suppose that 
 that $Q_1$ is not isomorphic to $\Bbb Z_2$.
 Then
 every normal elementary abelian subgroup 
 of dimension $r_2(S)$
  is contained in $Q_1\times Q_2$.
\end{enumerate}
}

\smallskip

\smallskip
\noindent {\em Proof.}
Let  $E\subseteq S$ be  a normal elementary abelian 
 2-group of dimension $nr_2(S)$.
If $E \subseteq Q_1 \times Q_2$, applying Lemma 2.2, we have
 $nr_2(S) = 2nr_2(Q_1)$.
 We shall therefore assume that 
 $E$ is not a subgroup of $ Q_1 \times Q_2$.
It follows that $E$ has an element of order 2
 of the form
 $q_1q_2 x$. Applying Lemma 3.1, $E$ contains $x$.
Note that $Q_1 \cap Q_2=1$,
 $[Q_1,Q_2]=1$ and that the conjugation action
 of $x$ interchanges $Q_1$ and $Q_2$.
 Hence  $q_1q_2 x$ is of order 2
 if and only if 
$$q_1xq_2x=1.\eqno(3.6)$$
 Since $E$ is abelian and $q_1q_2x, x \in E$,  we have
$xq_1q_2x = q_1q_2x^2 = q_1q_2.$
 Hence
$$q_1 = xq_2 x.\eqno(3.7)$$
 Applying (3.6) and (3.7), we have $q_1 = q_1^{-1}$.
 Hence the order of $q_1$ is either 1 or 2.
 Since this is true for all $q_1\in Q_1$
 (for any $q_1 \in Q_1$, let $q_2 = xq_1^{-1}x,$
 then equation (3.6) holds), we conclude
 that $Q_1$ is an elementary abelian 2-group.
 It is now clear that 
 $nr_2(S)  = 2nr_2(Q_1)$. 
\smallskip

\noindent (i) 
Suppose that $S$ possesses a normal elementary
 abelian subgroup $F$ 
such that $F$ is not a 
 subgroup of $Q_1\times Q_2$.
 It follows from the above that $Q_1$ is elementary 
 abelian.
A contradiction.

\smallskip
\noindent (ii) 
Suppose that $S$ possesses a normal elementary
 abelian subgroup $E$ of dimension $nr_2(S)$ 
such that $E$ is not a 
 subgroup of $Q_1\times Q_2$.
 It follows from the above that $Q_1$ is elementary 
 abelian. Set
$$E = A\times B,\mbox{ where } A = E \cap (Q_1\times
Q_2).\eqno(3.8)$$
 Applying the proof of Lemma 3.2
 (equation (3.5)), the rank
 of $Q_1$ is 1. It follows that $Q_1\cong \Bbb Z_2.$
 A contradiction. Hence $E \subseteq Q_1\times Q_2$.
\qed

\smallskip
\noindent {\bf Proposition 3.4.} {\em
 Let $T$ be a $2$-group and 
 let $w_n(T)$ be the $n$-th wreath product
 of $T$ and $\Bbb Z_2$ $(w_1(T) = T \wr \Bbb Z_2,
 \cdots ,
 w_n(T) = w_{n-1}(T)\wr \Bbb Z_2)$. Then
 $r_2(w_n(T)) = 2^n r_2(T)$ and
  $nr_2(w_n(T)) = 2^n nr_2(T)$.}

\smallskip
\noindent {\em Proof.} 
 Since 
$r_2(w_k(T))\ge 2$ and 
$w_k(T)$ is not elementary 
 abelian whenever $k \ge 1$,
 we may apply 
 the second part of Lemma 3.2 and 
 (i) of Lemma 3.3 to conclude that  
 $r_2(w_n(T)) = 2^{n-1} r_2(w_1(T))$ and
  $nr_2(w_n(T)) = 2^{n-1} nr_2(w_1(T))$.
 Applying the first part of Lemma 3.2,
 we have $r_2(w_1(T) = 2r_2(T)$.
 Applying the first part of Lemma 3.3,
 we have  $nr_2(w_1(T) = 2nr_2(T)$.
This completes the proof of the proposition. \qed

\subsection{
 Ranks and normal ranks of 
 Sylow 2-subgroups of $GL_n(q)$, $U_n(q)$, $O^+_{2n+1}
(q)$, $O_{2n}(\eta, q)$} Sylow 2-subgroups of 
 the above groups are determined by Carter
 and Fong [CF]. They are described as the direct
product of  certain  wreath products.
 Their rank and normal rank can be 
 determined by applying Proposition 3.4.
 The result of our study can be found in  the
 introduction.

\section{Twisted wreath product}

\noindent 
 {\bf Definition 4.1.} (Twisted wreath product)
Let $T\rtimes R$ be a semidirect product of 
 2-groups, where $T$ is normal.
 Define
 $$  w_1(T,R, J) =  \left < 
\left (\begin{array}
{cc}
T&0\\
0&1\\
\end{array}\right ),
\left (\begin{array}
{cc}
1&0\\
0&T\\
\end{array}\right ),
\left (\begin{array}
{cc}
r&0\\
0&r^{-1}\\
\end{array}\right ),
\left (\begin{array}
{cc}
0&1\\
1&0\\
\end{array}\right )\,:\, r\in R \right >.$$
 The above group 
is called the first twisted
 wreath
 product of $T$ and $R$ ($w_0(T,R, J) = T)$, where 
$$J =\left < 
\left (\begin{array}
{cc}
0&1\\
1&0\\
\end{array}\right )\right > \cong \Bbb Z_2.$$
Define similarly  the second
 twisted wreath product
 $w_2(T,R, J)$ to be the following group :
$$ \left < 
 \left (\begin{array}
{cc}
w_1(T,R,J)   &0\\
0&1\\
\end{array}\right ),
 \left (\begin{array}
{cc}
1&0\\
0&w_1(T,R,J) \\
\end{array}\right ),
\left (\begin{array}
{cc}
d_{2}(r)&0\\
0&d_{2} (r)^{-1}   \\
\end{array}\right ),
\left (\begin{array}
{cc}
0&I_2\\
I_2&0\\
\end{array}\right )\right >,$$
 where $d_{2}(r) $ is the $2\times 2$ diagonal
 matrix diag$\,(r ,1)$, $r \in R$.
One may define inductively $w_{n+1}(T,R,J)$ to be 
the following  (see [W] for more detail) :

\smallskip
{\small $$ \left < 
 \left (\begin{array}
{cc}
w_n(T,R,J)   &0\\
0&1\\
\end{array}\right ),
 \left (\begin{array}
{cc}
1&0\\
0&w_n(T,R,J) \\
\end{array}\right ),
\left (\begin{array}
{cc}
d_{n}(r)&0\\
0&d_{n}(r)^{-1} \\
\end{array}\right ),
\left (\begin{array}
{cc}
0&I_{2^n}\\
I_{2^n}&0\\
\end{array}\right )\right >,$$}
\smallskip

\noindent 
 where $d_{n}(r)$ is the $2^n\times 2^n$ diagonal
 matrix diag$\,(r,1,\cdots , 1)$,  $r \in R$.
 Define 
$$ R_{n+1} = \left <\mbox{diag}\,(d_n(r), d_n(r)^{-1})
 \right >\,,\,\,
 D_{n+1} = \left <
d_{n+1}(r)\,:\, r\in R\right >.
\eqno(4.1)$$

\subsection {Subgroups of 
 $w_{n+1}(T, R, J)$} 
 We shall define some subgroups
 (inductively) 
 of  $w_{n+1}(T, R, J)$ as follows :

\smallskip
\noindent 
(i) $ w_0((T, 1, J) = T$, and $w_{n+1}(T, 1, J)$
 is defined to be 

$$        \left < 
 \left (\begin{array}
{cc}
w_n(T,1,J)   &0\\
0&1\\
\end{array}\right ),
 \left (\begin{array}
{cc}
1&0\\
0&w_n(T,1,J) \\
\end{array}\right ),
\left (\begin{array}
{cc}
0&I_{2^n}\\
I_{2^n}&0\\
\end{array}\right )\right >.\eqno(4.2)$$

\smallskip
\noindent 
(ii) $w_0(T, R, 1) = T$, and 
 $w_{n+1}(T, R,1)$
 is defined to be 
 
$$  \left < 
 \left (\begin{array}
{cc}
w_n(T,R,1)   &0\\
0&1\\
\end{array}\right ),
 \left (\begin{array}
{cc}
1&0\\
0&w_n(T,R,1) \\
\end{array}\right ),
\left (\begin{array}
{cc}
d_{n}(r)&0\\
0&d_{n}(r)^{-1} \\
\end{array}\right )\right >.\eqno (4.3)$$

\smallskip
\noindent
 (iii) $w_{0}(1, 1, J) =1$, and 
$w_{n+1}(1,1,J)$
 is defined to be 
 
$$    \left < 
 \left (\begin{array}
{cc}
w_n(1,1,J)   &0\\
0&1\\
\end{array}\right ),
 \left (\begin{array}
{cc}
1&0\\
0&w_n(1,1,J) \\
\end{array}\right ),
\left (\begin{array}
{cc}
0&I_{2^n}\\
I_{2^n}&0\\
\end{array}\right )\right >.\eqno(4.4)$$

\smallskip
\noindent
 (iv) $w_{0}(T, 1, 1) =T$, and 
$w_{n+1}(T,1,1)$
 is defined to be 
 
$$    \left < 
 \left (\begin{array}
{cc}
w_n(T,1,1)   &0\\
0&1\\
\end{array}\right ),
 \left (\begin{array}
{cc}
1&0\\
0&w_n(T,1,1) \\
\end{array}\right ) \right >.\eqno(4.5)$$

\smallskip
\noindent It is clear that 
\begin{enumerate}
\item[(v)] $w_{n}(T, 1, J)$
 is just the $n$-th wreath product 
 of $T$ and $ J\cong \Bbb Z_2$. Equivalently,
$w_n(T,1,J) \cong w_n(T)$.
\item[(vi)] $w_{n}(T, R, 1)$
 is a normal subgroup of  $w_{n}(T, R, J)$.
\item[(vii)]
$w_{n} (1,1,J)$
consists of permutation matrices only. It is 
 isomorphic to a Sylow 2-subgroup 
 of $S_{2^{n}}$. 
\item[(viii)]
$w_{n} (T,1,1)$
 is a direct product of $2^{n}$ 
 copies of $T$'s and $w_{n} (T,1,1)$
 is a normal 
 subgroup of $w_{n} (T,R,J)$.

\end{enumerate}
It is also clear that 
$$w_n(T,R,J)
 \subseteq w_n(T,R,J)D_n = w_n(TR, 1,J).\eqno (4.6)$$

\subsection
 {Sylow 2-subgroups of the classical groups}
The main purpose of this section is the  construction
 of the group we define in (4.7).  We will
 see later (in section 5) that  Sylow 2-subgroups
 of the classical groups can be described by (4.7).
Let $n = 2^{m_1} + 2^{m_2} +\cdots +
2^{m_u}$ $(m_1 < m_2 <\cdots)$ 
be the 2-adic representation of $n$.
 In the case $u=1$, we define $U(R)=1$. In the 
 case $u \ge 2$, 
define $U(R)$ to be 
 the group generated by 
$$\mbox{diag}\,
 (d_{m_1}(r), d_{m_2}(r),I_{2^{m_3}},
I_{2^{m_4}},
\cdots, I_{2^{m_u}}   ),$$
$$\mbox{diag}\,
 (I_{2^{m_1}}, d_{m_{2} }(r), d_{m_{3}}(r),
 I_{2^{m_4}},\cdots , I_{2^{m_u}}   ),$$
$$\cdots \cdots\cdots \cdots\cdots
 \cdots\cdots \cdots \cdots \cdots\cdots\cdots \cdots$$
$$\cdots \cdots\cdots \cdots\cdots
 \cdots\cdots \cdots \cdots \cdots\cdots\cdots \cdots$$
$$\cdots \cdots\cdots \cdots\cdots
 \cdots\cdots \cdots \cdots \cdots\cdots\cdots \cdots$$
$$\mbox{diag}\,
 ( I_{2^{m_1}}   ,
\cdots , I_{2^{m_{u-2}}
},   d_{m_{u-1}}(r), d_{m_u}(r)),$$
 where $r\in R$. 
It is clear that $U(R)$ normalises $W(T,R,J)$,
 where $W(T,R,J)$ is defined to be 
 the following group.
$$ W(T,R,J) = \mbox{diag}\,
 (w_{m_1}(T,R,J),
 w_{m_2}(T,R,J),
\cdots , 
 w_{m_u}(T,R,J)).$$
 Define
$$S(T,R,J)  = W(T,R,J) \rtimes U(R).\eqno (4.7)$$
Note that we may define similarly the groups
 $W(T,R,1)$, $W(1,1,J), S(T,R,1)$ to be
$$W(T,R,1) = \mbox{diag}\,
 (w_{m_1}(T,R,1),
 w_{m_2}(T,R,1),
\cdots , 
 w_{m_u}(T,R,1)),$$
$$W(1,1,J) = \mbox{diag}\,
 (w_{m_1}(1,1,J),
 w_{m_2}(1,1,J),
\cdots , 
 w_{m_u}(1,1,J)),$$
and
$$S(T,R,1) =  W(T,R,1) \rtimes U(R)$$
respectively.
It follows from (4.6) that 
$S(T,R,J)\subseteq  S(TR,1,J)
 =    W(TR,1,J).$
 By (v) of the above, one has
$$W(TR,1,J) \cong \prod_{i=1}^u w_{m_i}(TR).
\eqno(4.8)$$

\noindent {\bf Remark.}   An easy observation of the 
 elements in $R_{n}$ (see (4.1)) and $U(R)$ shows that 
an element
 diag$\,(x_1, x_2, \cdots) \in W(TR,1,J)$
 is in $S(T,R,J)$ only if the numbers of the $i$'s
such that $x_i = t_i r_i$ ($t_i \in T$, $r_i
 \in R$, $r_i \ne 1$)  is even.

\subsection {Normal elementary abelian subgroups of
$S(T,R,J)$}
 Let $E$ be a normal elementary abelian subgroup of
$S(T,R,J)$ of dimension $nr_2(S(T,R,J))$
 The purpose of this section is to 
 give a supergroup of $E$ (Lemma 4.2) and a  lower
 bound of $r_2(E)$ (Lemma 4.3).
Proposition 4.4 is an easy consequence of 
 Lemmas 4.2 and 4.3 which works well
 for the Sylow 2-subgroups of the 
 classical groups.

\smallskip
\noindent 
 {\bf Lemma 4.2.} {\em Suppose that $T $ is not 
 elementary abelian.
Let $E$ be a normal elementary abelian subgroup
 of $S(T,R,J)$.
 Then $E \subseteq 
S(T,R,1) = 
 W(T,R,1) \rtimes U(R)$.
 In particular, 
 $r_2(E) \le n\cdot   r_2(T) + (n-1) r_2(R)$.}

\smallskip
\noindent {\em Proof.} Suppose not. Then
 $E$ has an element of the form 
 $w u \sigma$, where $w\in  W(T,R,1)$,
 $u \in U(R)$ and $\sigma \in W(1,1,J)-\{1\}$.
 Since the order of 
 $w u \sigma$ is 2, $\sigma$ is of order 2.
It is easy to see that 
 $\tau = w u \sigma$ normalises
 $S(T,1,1)$. Note that $S(T,1,1)$
 is the direct product of $n$ copies of $T$
 and that the conjugation of $\tau$
 permutes these $T$'s with a
 cycle decomposition $2^m$ for some $m\ge 1$.
 Without loss of generality, we may assume
 that $\tau$ interchanges the first two
 copies of the $T$'s. Label them as $Q_1$
 and $Q_2$, one has
 $$\tau \in  N = E \cap (Q_1\times Q_2)\rtimes
 \left <\tau\right >.$$ 
Since $N$ is a normal subgroup of 
$ (Q_1\times Q_2)\rtimes
 \left <\tau\right >$, we may apply (i) of Lemma 
 3.3 to  conclude that $N \subseteq Q_1\times Q_2$.
 A contradiction. Hence 
 $E \subseteq  W(T,R,1) \rtimes U(R)$.
Applying Lemma 2.3, it is clear that 
$$r_2(E) 
\le  r_2(W(T, R,1) \rtimes U(R))
\le  r_2(W(T, R,1)) + r_2(U(R)).$$
By Lemma 2.3,
$$ r_2(W(T, R,1))\le \sum_{i=1}^u
 \left( 2^{m_i} r_2(T) + (2^{m_i} -1) r_2(R)
\right)
\le  n\cdot r_2(T) + (n-u)r_2(R).
 $$ 
Note also that $r_2((U(R))= (u-1)r_2(R)$.
As a consequence,
$$r_2(E) \le  
   n\cdot  r_2(T) + (n-1) r_2(R).\eqno \qed $$

\smallskip
 
Let 
  $E_0$ and  $R_0$ be  
 elementary
 abelian subgroups of $T$ and $R$ respectively
 such that $\left <E_0, R_0\right>
 = E_0\times R_0$ is elementary abelian normal
 in $TR$.
 For $t\in T$, $r\in R_0$, 
Since 
 $\left < E_0, R_0\right>$
 is elementary abelian normal in $TR$,
 $t^{-1}rtr^{-1}  \in E_0\times R_0$. On the other 
 hand, since $R_0$ normalises $T$,
   $t^{-1}rtr^{-1} 
 \in T$. Hence 
 $t^{-1}rtr^{-1} \in (E_0\times R_0) \cap T = E_0$.
Hence
 $$
[T, R_0] \subseteq E_0.\eqno(4.9)$$
 Note also that
$R_0$ is normal in $R$ ($R$ is abelian) and that 
 the normality of $E_0\times 
 R_0$ in $TR$ implies that $E_0 \triangleleft T$.
 Hence 
$$ R_0 \triangleleft R\,,\,\,
    E_0 \triangleleft T.\eqno (4.10)$$

\smallskip
\noindent {\bf Lemma 4.3.} {\em  Let $E_0$ and
 $R_0$ be given as in the above and 
let $E$ be a normal elementary abelian subgroup
 of $S(T,R,J)$ of dimension
 $nr_2(S(T,R,J))$.
 Then 
 $r_2(E) \ge  n \cdot r_2(E_0) + (n-1) r_2(R_0)$.}

\smallskip
\noindent 
{\em Proof.} Let  
 $N = S(E_0, R_0,1)$. Since $[E_0, R_0]=1$,
 $N$ is elementary abelian.
Applying (4.9) and (4.10), one can show
 that
 $N$ is a normal subgroup of  $S(T,R,J)$.
$$r_2(N) 
= r_2(W(E_0, R_0,1) \times U(R_0))
 = r_2(W(E_0, R_0,1)) + r_2(U(R_0))$$
It is easy to see that 
$$ r_2(W(E_0, R_0,1))= \sum_{i=1}^u
 \left( 2^{m_i} r_2(E_0) + (2^{m_i} -1) r_2(R_0)
\right)
= n\cdot r_2(E_0) + (n-u)r_2(R_0).
 $$ 
Note that $r_2((U(R_0))= (u-1)r_2(R_0)$.
 This implies that 
$r_2(N) = n\cdot  r_2(E_0) + (n-1) r_2(R_0)$.
As a consequence,
$$r_2(E) \ge r_2(N)
 =  n\cdot   r_2(E_0) + (n-1) r_2(R_0).\eqno \qed $$

Applying Lemmas 4.2 and 4.3, we have the following
 proposition. Note that this proposition has its
 importance
as many Sylow 2-subgroups of  the classical groups 
 satisfy the assumptions of this proposition.

\smallskip
\noindent 
 {\bf Proposition 4.4.} 
 {\em 
Let $T$ and $R$ be given as in Definition $4.1$.
Suppose that $T$ is not elementary
 abelian. Then the following hold.
\begin{enumerate}
\item[(i)] Suppose that 
 there exists $E_0\triangleleft
  T$, $R_0 \triangleleft  E$
$(E_0$ and $R_0$ are elementary abelian with
 $\left <E_0,R_0\right >
 = E_0\times R_0)$
 such that   $r_2(T) = r_2(E_0),
\, r_2(R) = r_2(R_0)$. Then
 the normal rank of $S(T,R,J)$
 is $n\cdot r_2(T) + (n-1) r_2(R)$.
\item[(ii)] Suppose that 
 for every involution 
 $x \in TR-T $, there exists
 some $t \in T$ such that
 $\left < x, txt^{-1}\right >$
 is not abelian. 
 Let $E$ be a normal
 elementary abelian subgroup of 
 $S(T,R,J)$. Then 
 $E \subseteq S(T,1,1) = W(T,1,1)    $.
In particular, 
$
 r_2(E) =  
n \cdot nr_2(M),
$ where $M$ is the largest $($in dimension$)$
 elementary abelian normal subgroup of $TR$ that 
 contains in $T$.
\end{enumerate}
} 
\smallskip

\noindent {\em Proof.} (i) is an immediate
 consequence of Lemmas 4.2 and 4.3.
 Suppose that $TR$ satisfies the assumption of (ii).
 Then  normal
 elementary abelian subgroups of $TR$
 are subgroups of $T$.
By Lemma 4.2,
 $E \subseteq S(T,R,1) = W(T,R,1)\rtimes U(R)$
 (note that members in $S(T,R,1)$ are
 of the forms diag$\,(\cdots \cdots)$).
Let
 $$ u = \mbox{diag}\,
(e,e_2, \cdots , e_n) \in E.$$Suppose that 
 $e\notin T$. By our assumption, there exists $x\in
 T$ such that 
 $\left < e, txt^{-1}\right >$ is not abelian.
Let $$ v = \mbox{diag}\,
(x,x_2, \cdots , x_n) \in W(T,R,1)\rtimes U(R)
 \subseteq S(T,R,J).$$
 Then $\left <u, vuv^{-1}\right > \subseteq E$ is not abelian.
 A contradiction. Hence $e\in T$.
 Repeat this argument for
 all the $e_i$'s, we conclude that 
 $E \subseteq S(T,1,1)$.
It is now an easy matter to see that $
 r_2(E) = 
n \cdot nr_2(M)
.$
\qed

\subsection {Elementary abelian subgroups of
$S(T,R,J)$}
 Let $E$ be an elementary abelian subgroup of
$S(T,R,J)$ of dimension $r_2(S(T,R,J))$. 
 The purpose of this section is to 
 give upper and  lower
 bounds of $r_2(E)$ (Lemmas 4.5-4.8).
Proposition 4.9 is an easy consequence of 
 Lemmas 4.5-4.8 which works well
 for the Sylow 2-subgroups of the 
 classical groups.

 \smallskip
Let 
  $E_0$ and  $R_0$ be  
 elementary
 abelian subgroups of $T$ and $R$ respectively
 such that $\left <E_0, R_0\right>
 = E_0\times R_0$ is elementary abelian 
 in $TR$. Similar to Lemma 4.3, we have the following :

\smallskip
\noindent {\bf Lemma 4.5.} {\em 
Let $E$ be an elementary abelian subgroup
 of $S(T,R,J)$ of dimension $r_2(S(T,R,J))$.
 Then 
 $r_2(E) \ge  n\cdot  r_2(E_0) + (n-1) r_2(R_0)$.}

\smallskip
\noindent {\bf Lemma 4.6.} {\em
 $r_2(S(T,R,J))\le  
 n\cdot r_2(T) + (n-1) r_2(R).$}

\smallskip
\noindent
{\em Proof.} By Lemma 2.3,
$$r_2(S(T,R,J))
 \le r_2(W(T,R,J)) + r_2(U(R))
 =  \sum r_2(w_{m_i}(T,R,J))
 + (u-1)r_2(R).$$
By Lemma 2.3 and the definition of $w_{n+1}(T,R,J)$, we
have
 $r_2(w_{n+1}(T,R,J)
)\le 2r_2(w_n(T,R,J))+r_2(R) + 1$.
 Suppose that 
 $$r_2(w_{n+1}(T,R,J)
) =  2r_2(w_n(T,R,J))+r_2(R) + 1.\eqno (4.11)$$
Let $E$ be elementary abelian of dimension $r_2(S(T,R,J))$.
It follows from the equality (4.11) that 
 $E$ contains an element of the form
$\tau = g
 \left (\begin{array}
{cc}
0  &I_{2^n}\\
I_{2^n}& 0 \\
\end{array}\right ).$
Consider the conjugation action of $\tau$
 on 
$$ F = E \cap
 \left <
Q_1 = \left (\begin{array}
{cc}
w_n(T,R,J)  & 0\\
 0&  1 \\
\end{array}\right ),
Q_2 = \left (\begin{array}
{cc}
1  & 0\\
 0&  w_n(T,R,J) \\
\end{array}\right ) \right  >
.$$ 
Let $xy$ ($x\in Q_1$, $y\in Q_2$) be 
 an element of $F$. Since $F$ is abelian, one has
$$xy = \tau^{-1}x   y \tau = ( \tau^{-1}x 
\tau )(\tau^{-1}  y \tau).$$
 Since the action of $\tau $ interchanges
 $Q_1$ and $Q_2$ and $Q_1\cap Q_2 =1$, we have
$x = \tau^{-1}  y \tau.$
 As a consequence, $y$ is determined  uniquely
 by $x$. Hence the rank of $F$ is 
 at most $r_2(w_n(T,R,J))$.
On the other hand,  (4.11) and Lemma 2.3
 imply that $r_2(F) = 2 r_2(w_n(T,R,J)).$
 A contradiction. Hence 
(4.11) is false. It follows that 
 $$r_2(w_{n+1}(T,R,J)
) \le   2r_2(w_n(T,R,J))+r_2(R).\eqno (4.12)$$
 Hence 
 $$\sum r_2(w_{m_i}(T,R,J))
 \le \sum (2^{m_i} r_2(T) + (2^{m_i}-1) r_2(R)
= n\cdot r_2(T) + (n-u)r_2(R)
.$$
As a consequence,
 $$r_2(S(T,R,J))\le
 n\cdot r_2(T) + (n-1) r_2(R) .\eqno\qed$$

\smallskip
\noindent {\bf Lemma 4.7.} {\em
$r_2(S(T,R,J))\le  
 n\cdot r_2(TR).$
 Suppose that $r_2(T) = r_2(TR)$.
 Then $r_2(S(T,R,J))
 = n \cdot r_2(T)$.}

\smallskip
\noindent {\em Proof.}
Since $ w_n(TR, 1, J)$ is a wreath
 product
 (see (i) of section 4.1)), we may 
 apply Proposition 3.4 and conclude that 
 $$
  r_2(w_n(TR, 1, J)) =  2^n r_2(TR).
\eqno(4.13)$$
By (4.8),
 $$S(T,R,J) \subseteq S(TR, 1 ,J)
 = W(TR, 1,J).\eqno(4.14)$$
 By (4.13) and (4.14), we have
{ $$r_2(S(T,R,J))
 \le r_2( W(TR, 1 ,J))
= \sum 2^{m_i} r_2(TR)=
 n\cdot r_2(TR).$$}
It is clear that  
$$r_2(S(T,R,J)) \ge r_2(S(T,1,1)) = 
 n\cdot r_2(T).$$
This completes the proof of the lemma. \qed

\smallskip
\noindent {\bf Lemma 4.8.} {\em
 Suppose that $r_2(TR) = r_2(T) +1
 \ge 3$. Suppose further
 that there exists $E_0\subseteq 
  T$, $R_0 \subseteq   E$
 $(E_0$ and $R_0$ are elementary
 abelian$)$
 such that  $[E_0,R_0]=1$,
 $r_2(E_0) = r_2(T),$
 $r_2(R_0) =1$.
 Then $ r_2(S(T,R,J))= n\cdot r_2(TR) -1.$}

\smallskip
\noindent
{\em Proof.}
Since $E_0$, $R_0$ are elementary abelian
 with $[E_0, R_0]=1$, it follows that  
 $r_2(S(E_0, R_0,1)  )
   =n\cdot  r_2(E_0) + (n-1)r_2 (R_0)$.
 Since $r_2(TR) = r_2(T) +1$,   $r_2(E_0) = r_2(T),$
 $r_2(R_0) =1$, we have
 $$ r_2(S(T,R,J))\ge 
   n\cdot  r_2(E_0) + (n-1)r_2 (R_0)
 = n \cdot r_2(TR) -1.\eqno(4.15)$$

\smallskip
By Lemma 3.2, 
 $r_2(W(TR,1,J))
 = n\cdot  r_2(TR)$ and an 
 elementary abelian subgroup of 
$W(TR,1,J)$ of dimension $r_2(W(TR,1,J))$  is a subgroup of 
 $W(TR,1,1)$. As a consequence,
if $E$ is  elementary abelian subgroup of
dimension $r_2(W(TR,1,J))$,
 then
 $$E =\mbox{diag}\,
 (E_1, E_2, \cdots , E_{n}),$$
 where $r_2(E_i) = r_2(TR)$ and 
 the $E_i$'s are permuted by
$W(1,1, J)$. Since
 $r_2(TR) = r_2(T) +1 > r_2(T)$, $E_1$
 is not a subgroup of $T$. Hence 
 $E_1$ must admit an element of 
 the form $tr$, where $r \in R-\{1\}$.
Note that 
$$\mbox{diag}\,
 (tr, 1  \cdots ,  1) \in E$$
 is not in $S(T,R,J)$ (see remark of section 4.2).
 As a consequence, the rank
 of  $S(T,R, J)$ is less than $r_2(E) =
 n\cdot  r_2(TR)$.
Applying (4.15), we may now conclude that 
$
 r_2(S(T,R,J))= n\cdot r_2(TR) -1.$ \qed
 
\smallskip
The following is a consequence of Lemmas
 4.5-4.8. It is useful in the study
 of the ranks of Sylow 2-subgroups of 
the classical groups.

\smallskip
\noindent {\bf Proposition  4.9.} {\em
 The following hold.
\begin{enumerate}
\item[(i)]
 Suppose 
 that there exists $E_0\subseteq 
  T$, $R_0 \subseteq   R$
 $(E_0$ and $R_0$ are elementary
 abelian$)$
 such that  $[E_0,R_0]=1$,
 $r_2(E_0) = r_2(T),$
 $r_2(R_0) =r_2(R)$. Then
 $$r_2(S(T,R,J)) = 
 n\cdot r_2(T) + (n-1) r_2(R).$$
\item[(ii)] Suppose that $r_2(TR) = r_2(T)$.
 Then $r_2(S(T,R,J)) = n\cdot r_2(T)$.
\item[(iii)]
 Suppose that $r_2(TR) = r_2(T) +1
 \ge 3$. Suppose further
 that there exists $E_0\subseteq 
  T$, $R_0 \subseteq   R$
 $(E_0$ and $R_0$ are elementary
 abelian$)$
 such that  $[E_0,R_0]=1$,
 $r_2(E_0) = r_2(T),$
 $r_2(R_0) =1$.
 Then $ r_2(S(T,R,J))= n\cdot r_2(TR) -1.$
\end{enumerate}
}

\section{Applications :  Ranks and Normal 
 ranks of Groups of Lie type}
\smallskip

\subsection{Special linear and special
 unitary groups I}
 Let $2^{t+1}$ be the largest power of 2 that
 divides $q^2-1$ and 
let $T$ be a generalised quaternion group
 $\left <v, w\right >$,
 where $o(v)= 2^{t}, o(w) = 4$.
Suppose that $T$ is normalised by $R = 
\left
 <e\right >$ (a 
 group of order 2) :
$$eve^{-1} = v^{- 1}\,, \,\,
ewe^{-1} = vw.$$
Let 
  $n = 2^{m_1} + 2^{m_2} + \cdots  +2^{m_u}$ 
 be the $2$-adic representation of $n$.
By Theorem 4 of [W], 
 $S(T,R,J)$ is a Sylow 2-subgroup of 
  $SL(2n,q)$
 $(q\equiv 3$ (mod 4)) and
 $SU(2n,q)$
 $(q\equiv 1$ (mod 4)).

\smallskip
\noindent  (a) {\bf   Normal Rank of  $S(T,R,J)$.} 
 One can show easily that 
 $x\in TR-T$ is an involution
 if $x$ is of the form $v^i e$. Further,
        $\left < x, wxw^{-1}\right >$
 is not abelian. Applying (ii) of Proposition  4.4,
 $$nr_2(S(T,R,J))
 = n\cdot r_2(T) = n.$$

\smallskip
\noindent (b)  {\bf    Rank of  $S(T,R,J)$.} 
One sees easily 
 that  $r_2(TR) = 2$, $r_2(T) = r_2(R) = 1$.
 Further, $E_0   =
\left <v^{2^{t-1}}\right >$ and $R_0 =\left < e\right >$ 
  satisfy the assumption of (i)
 of Proposition 4.9.
 Hence 
 $$r_2(S(T,R,J))  = 2n-1.$$

\smallskip
\noindent {\bf Remark.}
In the case $S$ is a Sylow 2-subgroup of 
 $SL(2n+1, q) $ ($q \equiv 3$ (mod
 4) or 
 $SU(2n+1, q) $ ($q \equiv 1$ (mod 4)
 it is well known that $S$ is 
isomorphic to  a Sylow 2-subgroup of $GL(2n,q)$. 
 Its rank and normal rank are given 
 in section 3.1.

\smallskip
\subsection{Special linear and special
 unitary groups II}
 Let $2^{t+1}$ be the largest power of 2 that
 divides $q^2-1$ and 
let $T$ be a generalised quaternion group
 $\left <v, w\right >$,
 where $o(v)= 2^{t}, o(w) = 4$.
Suppose that $T$ is normalised by $e$ (an 
 element of order $2^t$) :
$$eve^{-1} = v\,, \,\,
ewe^{-1} = vw.$$
Let 
  $n = 2^{m_1} + 2^{m_2} + \cdots  +2^{m_u}$ 
 be the $2$-adic representation of $n$.
Define 
 $  R=\left <e \right >$. By Theorem 4 of [W], 
$S(T,R,J)$ is
 a Sylow 2-subgroup of $SL(2n,q)$
 $(q\equiv 1$ (mod 4)) and
 $SU(2n,q)$
 $(q\equiv 3$ (mod 4)).
 Let  $E_0 = \left < w^2 = v^{2^{t-1}}, \right
 > \triangleleft T$, $R_0 =\left <f = e^{2^{t-1}}\right >
\triangleleft R$. Then $[E_0, R_0]=1$ and $E_0\times
 R_0$ 
 is normal elementary abelian in $TR$.

\smallskip
\noindent  (a) {\bf   Normal Rank of  $S(T,R,J)$.} 
 Since $r_2(T) = r_2(E_0)$, $r_2(R) = r_2(R_0)$,
 one may apply (i) of Proposition 4.4  and conclude
 that 
$$nr_2(S(T,R,J)) = 2n-1.$$

\smallskip
\noindent  (b) {\bf  Rank of  $S(T,R,J)$.} 
 By (i) of Proposition 4.9, we have 
$$r_2(S(T,R,J)) = 2n-1.$$

\smallskip
\noindent {\bf Remark.}
In the case $S$ is a Sylow 2-subgroup of 
 $SL(2n+1, q) $ ($q \equiv 1$ (mod
 4) or 
 $SU(2n+1, q) $ ($q \equiv 3$ (mod 4)
 it is well known that $S$ is 
isomorphic to  a Sylow 2-subgroup of $GL(2n,q)$. 
 Its rank and normal rank are given 
 in section 3.1.

\subsection{Symplectic groups}
 Let $T$ be a generalised quaternion group whose
 order is the largest power of 2 that divides
 $q^2-1$ and let $w_m$ be the wreath product of 
 $m$ copies of $\Bbb Z_2$. Then $T\wr w_m$ is 
 a Sylow 2-subgroup of $Sp_{2^m}(q)$.
 Let $n = 2^{m_1} +  2^{m_2}+ \cdots + 2^{m_u}$
 be the 2-adic representation of $n$ and
 let
 $T_i$ be a Sylow 2-subgroup of $Sp_{2n_i}(q)$
 ($n_i= 2^{m_i}$).  By Theorem 6 of [W], 
 $S = T_1 \times T_2 \times \cdots \times T_u$ is 
 a Sylow 2-subgroup of $Sp_{2n}(q)$.
We may 
 apply Proposition 3.4 to conclude that 
 the rank as well as the normal
 rank of $S$ 
is  $n$.

\subsection{Orthogonal Commutator Groups $\Omega
 _{2n+1}( q)$} Let $2^{t+1}$ be the largest power
 of 2 that divides $q^2-1$ and let $T =\left < v,w\right>
$ be a dihedral group of order $2^{t}$,
 where $o(v)= 2^{t-1}, o(w)= 2, wvw= v^{-1}$.
 Further, $R = \left < e\right >$ is a group of 
 order 2 acts on $T$ by
 $o(e)=2$, $eve= v^{-1}, ewe= vw.$
Let 
  $n = 2^{m_1} + 2^{m_2} + \cdots  +2^{m_u}$ 
 be the $2$-adic representation of $n$.
 Then $S(T,R,J)$  is 
 a Sylow 2-subgroup of $\Omega_{2n+1}(q)$
 (see (ii) of Theorem 7 of Wong [W]).

\smallskip

\noindent {\bf Case 1. $T$ is nonabelian.}

\begin{enumerate}
\item[(a)]
 {\bf  Normal rank of  $S(T,R,J)$.} 
Direct calculation shows that $de$ $(d\in T, 
 e\in R$) is an involution
 if and only if $ede=d^{-1}$. As a consequence,
 $d= v^i$ for some $i$.
 Hence involutions in $TR$ take the forms $v^ie$.
 One sees easily 
 that if $x= v^ie \in TR-T$, then  
 $\left <x, wxw^{-1}\right >$ is not abelian.
By  (ii) of Proposition 4.4,
 $$nr_2(S(T,R,J)
 = n\cdot nr_2(TR) = n.$$

\noindent {\bf Remark.} 
In the case $|T| \ge 16$, the normal rank
 of $T$ is 1. It follows that $nr_2(S) = n$.
 In the case $|T| = 8$, the normal rank of 
 $T$ is 2. However, such normal subgroups
 can not be normalised by $e$.
 Hence $nr_2(TR) = 1$.

\item[(b)]   {\bf  Rank of  $S(T,R,J)$.} 
 We shall now determine $r_2(S)$. 
Recall that $de \in TR$ is of order 2 only if 
 $d= v^i$. It follows that 
 $r_2(TR) = r_2(T) =2$. We may apply
 (ii) of Proposition 4.9  and 
 conclude that
 $$r_2(S(T,R,J)) = 2n.$$
\end{enumerate}

\smallskip
\noindent {\bf Case 2. $T$ is abelian.}
 This implies that $T \cong  E_4 = \Bbb Z_2 \times \Bbb Z_2$.
Consequently, we have
 $o(v)=o(w)=o(e) = 2$, $vw=wv$, $ev=ve$, $ewe=vw$.
It follows that $S(T,R,J)$ is isomorphic to 
 a Sylow 2-subgroup of the alternating group
 $A_{4n}$. It is well known that 
 the rank of $A_{4n}$ is $2n$
 (see Proposition 5.2.10 of 
 [GLS3]).
 Hence $S(T,R,J) $ is of rank $2n$ as well.
Note that $S(T,1,1)$ is elementary abelian
 normal of rank $2n$. Hence 
  the normal rank of
 $S(T,R,J)$ is $2n$.

\subsection{Orthogonal Commutator Groups $\Omega
 _{2n}(\eta, q)$ where 
 $\eta = \pm 1$, $q^n \equiv -\eta$
 (mod 4)}
Applying Theorem 7 of Wong
 [W],  a Sylow 2-subgroup of 
 $\Omega_{2n}(\eta, q)
 = P\Omega_{2n}(\eta, q)$ is isomorphic 
 to a Sylow 2-subgroup of 
 $O_{2(n-1)} (\eta ',  q)$, where
 $q^{n-1}\equiv \eta '$ (mod 4). Let
 $S$ be a Sylow 2-subgroup of  $O_{2(n-1)} (\eta ',
 q)$, where
 $q^{n-1}\equiv \eta '$ (mod 4).
Applying Theorem 3 of Carter and Fong [CF],
 $S$ is isomorphic to a Sylow
 2-subgroup of $O^+_{2n-1}(q)$.
We shall now describe $S$  as follows : 

\smallskip
\noindent 
Let $D$ be a dihedral group of order $2^{s+1}$,
 where $2^{s+1}$ is the largest power of 2
 that divides $q^2-1$. Then 
 $D$ is isomorphic to a Sylow 2-subgroup of 
 $O^+_3(q)$. Let $T_{r-1}$ be the wreath product
 of $r-1$ copies of $\Bbb Z_2$ and let
 $S_r$ be the wreath product of $D$ and $T_{r-1}$.
 Then $S_r$ is a Sylow 2-subgroup of 
 $O_{2^r+1}^+(q)$.
 Let $2(n-1)= 2^{m_1} +  2^{m_2}+ \cdots + 2^{m_u}$
be the 2-adic representation of $2(n-1)$. Applying 
 Theorem 2 of Carter and Fong [CF],
 $$S \cong
 S_{m_1} \times  S_{m_2}\times \cdots \times  S_{m_u}.$$
Note that the normal rank of $D$ is 1 if $D$ is 
 of order 16 or more and 2 if the order of $D$ is 8.
The normal rank of 
 $S$ can be determined by applying Proposition 3.4.
In particular, one has
\begin{enumerate}
\item[(i)]
 if $|D| \ge 16$ , then $nr_2(S) = n-1$,
 \item[(i)]
 if $|D|  = 8$ , then $nr_2(S) = 2n-2$
\end{enumerate}
 
\smallskip
\noindent 
By Proposition 3.4, the rank of $S$ is $2n-2$.

\subsection{Orthogonal Commutator Groups $\Omega
 _{2n}(\eta, q)$, where $n$ is even,
 $\eta = \pm 1$, $q^n \equiv \eta$
 (mod 4)} Note that $\eta =1$ since $n$ is even and 
 $q^n \equiv \eta$.
 Let $2^{t+1}$ be the greatest power
 of 2 that divides $q^2-1$ and let
$T$ be the central product of two dihedral
 groups of order $2^{t+1}$ :
$$T = \left
< d= \left (\begin{array}
{cc}
u&0\\
0&u^{-1} \\
\end{array}\right ),
g=\left (\begin{array}
{cc}
u&0\\
0&u \\
\end{array}\right ),
h= \left (\begin{array}
{cc}
0&1\\
1&0 \\
\end{array}\right ),
k=\left (\begin{array}
{cc}
0&w\\
w&0 \\
\end{array}\right )\right>,$$
 where
$o(u) = 2^t$, $o(w)=2$, $wuw=u^{-1}$.
Let
$ R = \left < e, f\right > \cong \Bbb Z_2
 \times \Bbb Z_2$, where 
$$d^e=g^{-1}, g^e=d^{-1}, h^e=gk, k^e=dh,$$
$$d^f=g, g^f= d, h^f=k,k^f=h.$$
Let 
 $n/2 = 2^{m_1} +  2^{m_2}+\cdots +  2^{m_u}$ $(
 n$ even) be 
 the 2-adic representation of $n/2$.
By Theorem 11 of [W], 
 $S(T,R,J)$
  is a Sylow 2-subgroup of $\Omega_{2n}(\eta,
 q)$, where $n$ is even,
 $\eta = \pm 1$, $q^n \equiv \eta$
 (mod 4). 

\smallskip
\noindent  (a) {\bf  Rank of  $S(T,R,J)$.} 
One sees easily that 
 $r_2(TR) = 4 = r_2(T) + 1$.
Let  $E_0 = \left <dg, d^{-1}g, hk\right >$,
 $R_0 = \left < f\right >$. Then
 $r_2(E_0) = 3$, 
 $[E_0,R_0]=1$.
 Further, $S(E_0, R_0,1)$ is of rank $2n-1$.
Applying  (iii) of Proposition 4.9,
$$r_2(S(T,R,J))
 = (n/2)r_2(TR)  -1 = (n/2) 4-1 = 2n-1
.$$

\smallskip
\noindent {\bf Remark.} Since $[E_0,R_0]=1$,
$[S(E_0, R_0,1), \mbox{diag}\,(f, 1, \cdots)
 ] =1$.

\medskip
\noindent  (b) {\bf  Normal Rank of  $S(T,R,J)$.} 
 We shall now study the normal rank of 
$S(T,R,J)$.

\smallskip
\begin{enumerate}
\item[b(i)]
 {\bf Case 1. $t\ge 3$.} 
 Direct calculation
 shows that $  v = d^xh^yg^zk^we$ is of order 2 if and only if 
\begin{enumerate}
\item[(i)] $y\equiv w \equiv 0$ (mod 2),
 $x \equiv z$ (mod $2^{t-1})$,
 or
\item[(ii)] $y\equiv w \equiv 1$ (mod 2), $x+z-1 \equiv 0$
 (mod $2^{t-1})$.
\end{enumerate}
Further,$\left < v, dvd^{-1}\right >$ is not
 abelian.

\smallskip\noindent
In the case
$  r = d^xh^yg^zk^wf$, $r$  is of order 2 if and only if 
\begin{enumerate}
\item[(i)] $y\equiv w \equiv 0$ (mod 2),
 $x +z\equiv 0$ (mod $2^{t-1})$,
 or
\item[(ii)] $y\equiv w \equiv 1$ (mod 2), $x-z \equiv 0$
 (mod $2^{t-1})$.
\end{enumerate}
Further,$\left < r, drd^{-1}\right >$ is not
 abelian.

\smallskip\noindent
In the case
$  s = d^xh^yg^zk^wef$, $s$  is of order 2 if and only if 
 $y\equiv w \equiv 0$ (mod 2).
Further,$\left < s, dsd^{-1}\right >$ is not
 abelian.
\smallskip

\noindent By  (ii) of Proposition  4.4,
  $nr_2(S(T,R,J)) = n \cdot nr_2(TR)$.
It is easy
 to see that $ nr_2(TR) = 2$ and that 
$$E_0 =\left < (dg)^{2^{t-2}}, (d^{-1}g)^{2^{t-2}}
\right > \triangleleft  TR$$ is of rank 2.
Hence  $nr_2(S(T,R,J)) = (n/2)  nr_2(TR) =n$.
Further, $S(E_0,1,1)$ is elementary abelian normal of 
 dimension $n$.
\smallskip

\noindent {\bf Remark.} Note that our results about
 elements of order 2 in $TR-T$ works for $t=2$ as well.

\smallskip
\item[b(ii)] {\bf Case 2. $t = 2$.} In this case, $T$ is 
 a central product of two dihedral groups of order 8.
Note first that  $nr_2(T) = 3$.
Let   $E_0 =\left < dg, d^{-1}g, hk
\right > \triangleleft T$,
 $R_0=\left < f\right > \triangleleft R.$
Then $[E_0, R_0]=1$ and $E_0R_0$ 
 is normal elementary abelian of 
 rank 3 in $TR$. It follows that
$S(E_0, R_0, 1) $ is normal elementary
 abelian in $S(T,R,J)$ and that 
$$
r_2(S(E_0, R_0, 1))
 = (n/2) r_2(E_0) + (n/2 -1) r_2(R_0)
 = 2n-1.$$
 As a consequence,
 $nr_2(S) = 2n-1$ ((i) of Proposition 4.4).

\smallskip
\noindent {\bf Remark.} Since $[E_0, R_0]=1$,
 we have $[S(E_0, R_0,1), \mbox{diag}\,
 (f, 1, \cdots )] =1$. 

\end{enumerate}

\subsection{Orthogonal Commutator Groups $\Omega
 _{2n}(\eta, q)$ where $n$ is odd,
 $\eta = \pm 1$, $q^n \equiv \eta$
 (mod 4)}
 Let $n = 1+n_1$, where $n_1$ is even.
Let $S$ be a Sylow 2-subgroup of 
$\Omega
 _{2n_1}(\eta, q)$ where $n_1$ is even,
 $\eta = \pm 1$, $q^{n_1} \equiv \eta$
 (mod 4) and let 
 diag$\,(e,1,1,\cdots , 1) = d(e)$,
 diag$\,(f,1,1,\cdots , 1)= d(f)$
 ($e$ and $f$ are given in section 5.6).
Set
$$\left < S, d(e), d(f)\right >
 \times \left < x, y\right >,
\mbox{ where } o(x)=o(y) = 2, o(xy) = 2^t.$$
Let  $ D = \left < d(e)x,d(f)y\right > $
 (dihedral of order $2^{t+1}$).
 By Theorem 12 of Wong [W], 
 $$ V =  \left < S, d(e)x,d(f)y\right >
 = S \rtimes D \eqno(5.1)$$
 is a Sylow 2-subgroup of 
$\Omega
 _{2n}(\eta, q)$.

\smallskip
\noindent  (a) {\bf   Normal Rank of  $V$.} 
 Let $E$ be normal elementary abelian of dimension $nr_2(V)$. 
 Recall  that $ D = \left < d(e)x,d(f)y\right > $
 is dihedral of order $2^{t+1}$ $(t\ge 2)$.
Set
 $$E = (E\cap S) \times B.\eqno (5.2)$$
 By Lemma 2.3, $r_2(B) \le 2$.
 Hence $$
nr_2(V)=
r_2(E) \le  nr_2(S) + 2 .\eqno(5.3)$$
In the case $D$ is of order 8 ($t =2$), $D$ has a
 normal subgroup $
 \left < d(f)y, (xy)^2\right > =D_0 \cong \Bbb Z_2
\times
 \Bbb Z_2$. 
 Let  $K 
= S(E_0,R_0,1) 
\subseteq S$ be given as in b(ii) of section 5.6.
Then $K$ is 
 elementary abelian of dimension $nr_2(S)$.
 Since $\left <x,y
\right >$ commutes with $S$ and $[K, d(f)]=1$
 (see remark of b(ii) of section 5.6),
 $[K, D_0] =1$. Hence $K\times D_0$  is 
 normal
 elementary abelian (of $V$) of dimension  $nr_2(S) + 2$.
 Applying inequality (5.3), we have 
 $$nr_2(V) =nr_2(S) + 2
 = (2n_1-1)+2 = 2n_1+1 = 2n-1.$$
In the case, $|D| >8$ $(t\ge 3$). If $E$ is of 
dimension 
 $nr_2(S) +2$. Then $B$ (see (5.2)) is of dimension  2.
 Consequently, $B|D$ 
 (see Lemma 2.2 for the notation $B|D$) is elementary
 abelian of  dimension 2.
 Since
 $E$ is normal in  $ V
 = S\rtimes D$, one can show that 
 $E|D = B|D$ (of dimension 2) is normal in $D$. This is a
 contradiction (the normal rank of $D$ is 1).
 Hence (5.3) can be improved to 
$$nr_2(V) \le nr_2(S) + 1
 =n_1+1.\eqno (5.4)$$
Let  $K  =S(E_0,1,1) \subseteq S$ be  normal
 elementary abelian (of $S$) of dimension $nr_2(S)$
 (see b(i) of section 5.6) and let
 $D_0 = \left < (xy)^{o(xy)/2} \right >
 = Z(D)$.
  It is 
 clear that 
 $K\times D_0$ 
 ($K \subseteq S$ commutes with $xy$) is 
 normal
 elementary abelian (of $V$) of dimension  $r_2(S) + 1$.
 Applying inequality (5.4), we have 
 $$nr_2(V) =nr_2(S) + 1
 = n_1+1 =n.$$
\smallskip
In summary, we have
\begin{enumerate}
\item[(i)] if $t\ge 3$, then 
$nr_2(V) = n_1+1 = n$,
\item[(ii)] if $t = 2$, then 
$nr_2(V ) =  2n_1+1= 2n-1$.
\end{enumerate}

\smallskip
\noindent  (b) {\bf   Rank of  $V$.} 
We shall now determine the rank of 
$ V =  S\rtimes D$.
Let $E$ be elementary abelian 
 of $\left < S, d(e)x,d(f)y\right >$
 of dimension $r_2(V)$. 
 Let $$E = A\times B,\mbox{ where } A =
 E \cap S.$$
 Since $r_2(B) \le r_2(D)$ (Lemma 2.3)
 and $ r_2(S) = 2n_1-1$ (see (a) of section 5.6), we have
$$ r_2(V)= 
r_2(E) = r_2(A) + r_2(B) \le
 r_2(S) +2 = 2n_1+1.$$
 Let $ \left < d(f)y, 
 (xy)^{o(xy)/2} \right >
 = D_0 \subseteq D$. Then $D_0$ is   elementary
 abelian
 of dimension 2. Let $K =S(E_0,R_0,1)  \subseteq S$
 be elementary abelian of dimension  $r_2(S)
 = 2n_1-1$ (see (a) of section 5.6).
Applying 
the remark of (a) of section 5.6 and the fact that 
$[S, \left <x,y\right >]=1$,  $KD_0 = K\times D_0$ 
 is elementary abelian of dimension $2n_1+1$.
 Hence $r_2(V) = 2n_1+1 = 2n-1$.

\bigskip
\bigskip
{\small
\noindent DEPARTMENT OF MATHEMATICS,\\
NATIONAL UNIVERSITY OF SINGAPORE,\\
SINGAPORE 117543,\\
REPUBLIC OF SINGAPORE.}

\bigskip

\noindent {\tt e-mail: matlml@math.nus.edu.sg}

\bigskip
\bigskip

\noindent  lang-67-6.tex

\end{document}